\documentclass[12pt]{amsart} % use larger type; default would be 10pt

\usepackage[utf8]{inputenc} % set input encoding (not needed with XeLaTeX)
\usepackage{amsmath}
\usepackage{youngtab}
\usepackage{mathrsfs}
\usepackage{amssymb}

\usepackage{geometry} % to change the page dimensions
\usepackage{graphicx} % support the \includegraphics command and options

\usepackage{booktabs} % for much better looking tables
\usepackage{array} % for better arrays (eg matrices) in maths
\usepackage{verbatim} % adds environment for commenting out blocks of text & for better verbatim
\usepackage{subfig} % make it possible to include more than one captioned figure/table in a single float
\usepackage{chngcntr}

\usepackage{fancyhdr} % This should be set AFTER setting up the page geometry
\newtheorem{theorem}{Theorem}

\newtheorem{lemma}{Lemma}
\counterwithin{equation}{section}
\counterwithin{lemma}{section}
\counterwithin{theorem}{section}
\counterwithin{defn}{section}

\title{On Tensor Product Decomposition of $\widehat{\mathfrak{sl}}(n)$ Modules}
\author{Kailash C. Misra \and  Evan A. Wilson}
\address{Department of Mathematics, North Carolina State University,  Raleigh,  
NC 27695-8205}
\email{misra@ncsu.edu ,   eawilso4@ncsu.edu}
\thanks{Partially supported by NSA grant, H98230-08-1-0080.} 
\date{} % Activate to display a given date or no date (if empty),
\begin{document}
\maketitle
\begin{abstract}
We decompose the $\widehat{\mathfrak{sl}}(n)$-module $V(\Lambda_0) \otimes V(\Lambda_0)$ and give generating function identities for the outer multiplicities. In the process we discover some seemingly new partition identities in the cases $n=2,3$.
\end{abstract}
\section{Introduction}
The affine Lie algebra $\widehat{\mathfrak{sl}}(n)$ is the infinite dimensional analog of the finite dimensional simple Lie algebra $\mathfrak{sl}(n)$ of $n \times n$ trace zero matrices (cf. \cite{K}). For a dominant integral weight $\lambda$, let $V(\lambda)$ denote the irreducible integrable $\widehat{\mathfrak{sl}}(n)$-module with highest weight $\lambda$. The level of $V(\lambda)$ is the scalar $k\ge 0$ by which the canonical central element $c \in \widehat{\mathfrak{sl}}(n)$ acts. Around 1990, Kashiwara (\cite{Ka1}, \cite{Ka2}) and Lusztig \cite{Lu} independently introduced the crystal $\mathcal{B}(\lambda)$ associated with $V(\lambda)$. In \cite{JMMO} an explicit realization of the crystal $\mathcal{B}(\lambda)$ is given in terms of extended Young diagrams (or colored Young diagrams). Let $P_+$ denote the set of dominant integral weights. For $\lambda , \mu \in P_+$, it is known that the tensor product module $V(\lambda)\otimes V(\mu)$ is completely reducible (cf. \cite[Corollary 10.7 b]{K}), that is:

\begin{equation}
V(\lambda) \otimes V(\mu)=\bigoplus_{\nu \in P_+}V(\nu)^{\oplus b_{\nu}} \label{gendecomp},
\end{equation}
where $b_{\nu}$ denotes the number of times $V(\nu)$ occurs in this decomposition and is called its outer multiplicity. Although there are Littlewood-Richardson type combinatorial algorithms known (for example see \cite{Li}, \cite{OSS}) to decompose such tensor products, explicit determination of the outer multiplicities are known only for the case $n =2$ \cite{F}. In this paper we use the crystal consisting of $n$-regular colored Young diagrams for the integrable module $V(\Lambda_0)$ \cite{MM} and determine the irreducible summands in \eqref{gendecomp} for $\lambda = \Lambda_0 = \mu$ and their outer multiplicities. Furthermore, we use the principally specialized characters to determine certain generating functions of these outer multiplicities. Finally comparing these two approaches of computing the generating functions for the outer multiplicities in the cases $n =2, 3$, we obtain several combinatorial identities. Some of these identities seem to be new. 

\section{Preliminaries}

Let $E_i:=E_{i,i+1}$,  $F_i:=E_{i+1,i}$, and $ H_i:=E_{ii}-E_{i+1,i+1}$, $ i=1\dots n-1$ be the Chevalley generators of the finite dimensional simple Lie algebra $\mathfrak{sl}(n,\mathbb{C})$ of $n \times n$ trace zero matrices, where $E_{ij}$ denotes the matrix with a $1$ in the $i^{\text{th}}$ row and $j^{\text{th}}$ column, and 0 elsewhere.  Also, let $E_0:=E_{n1}$, $F_0:=E_{1n}$, and $H_0:=[E_0,F_0].$  The affine Lie algebra $\widehat{\mathfrak{sl}}(n) = \mathfrak{sl}(n,\mathbb{C}) \otimes \mathbb{C}[t,t^{-1}] \oplus \mathbb{C} c \oplus \mathbb{C} d$ is generated by the degree derivation $d=1 \otimes t\frac{d}{dt}$, and the Chevalley generators: 
$$ e_i=E_i \otimes 1,  \qquad f_i=F_i \otimes 1, \qquad  h_i=H_i \otimes 1,\qquad i=1,2,\dots, n-1, $$
$$ e_0=E_0 \otimes t,\qquad f_0=F_0 \otimes t^{-1}, \qquad h_0=H_0 \otimes 1 + c,$$
where $c=\sum_{i=0}^{n-1} h_i$ spans the one-dimensional center. 
The \emph{Cartan subalgebra} of  $\widehat{\mathfrak{sl}}(n,\mathbb{C})$ is the abelian subalgebra ${\mathfrak{h}}:=\text{span}_{\mathbb{C}}\{h_0, h_1, \dots, h_{n-1}\} \oplus \mathbb{C}d$.  Let $\Delta =\{\alpha_0, \alpha_1, \dots, \alpha_{n-1}\}$ be the set of simple roots and
$\Phi$ be the set of roots for  $\widehat{\mathfrak{sl}}(n,\mathbb{C})$. Then $\delta = \alpha_0+ \alpha_1+ \cdots + \alpha_{n-1}$ is the null root. 
The \emph{weight lattice} of $\widehat{\mathfrak{sl}}(n,\mathbb{C})$ is defined to be $P:=\text{span}_{\mathbb{Z}}\{\Lambda_0,\Lambda_1,\dots, \Lambda_{n-1}, \delta\},$ where $\langle \Lambda_i, h_j\rangle=\delta_{ij},$ and $\langle \Lambda_i,d\rangle=0$, for $i=0,1,\dots,n-1$.  The  \emph{dominant weight lattice} is defined as $P_+=\text{span}_{\mathbb{Z}_{\geq 0}}\{\Lambda_0,\Lambda_1,\dots,\Lambda_{n-1}\}\oplus \mathbb{Z}\delta.$  By the \emph{level} of $\lambda \in P_+$ we mean the nonnegative integer $\text{level}(\lambda)= \lambda (c)$.  For notational convenience we define $\alpha_i = \alpha_{\overline{i}}$ and $\Lambda_i = \Lambda_{\overline{i}}$ for all $i \in \mathbb{Z}$ where  $\overline{i}:=i \pmod{n}$.

Let $V(\Lambda_0)$ be the integrable highest weight module with highest weight $\Lambda_0$. Then $V(\Lambda_0) \otimes 
V(\Lambda_0)$ is \emph{completely reducible} (cf. \cite[Corollary 10.7 b]{K}), that is: 
\begin{equation}
V(\Lambda_0) \otimes V(\Lambda_0) = \bigoplus_{\lambda \in P_+}V(\lambda)^{\oplus b_{\lambda}} \label{decomp},
\end{equation}
where $b_{\lambda}\in \mathbb{Z}_{\geq 0}$ is called the \emph{outer multiplicity} of $V(\lambda)$ in $V(\Lambda_0) \otimes V(\Lambda_0)$.  In this paper, we give the highest weights of each of the irreducible summands and determine their outer multiplicities. 

%the outer multiplicities that occur.

The theory of crystals (cf. \cite{HK}) associated  with integrable representations of quantum affine algebras has provided an useful tool to study combinatorial properties. Indeed in this paper we will use the explicit realization of the crystal 
$\mathcal{B} = \mathcal{B}(\Lambda_0)$ for the module $V(\Lambda_0)$ in terms of \emph{extended Young diagrams} (or \emph{colored Young diagrams})(\cite{MM}, \cite{JMMO}) which we briefly describe.

Let $I:=\{0,1,\dots n-1\}$ denote the index set for $\widehat{\mathfrak{sl}}(n)$.  An \emph{extended Young diagram} is a collection of $I$-colored boxes arranged in rows and columns, such that the number of boxes in each row is greater than or equal to the number of boxes in the row below.  To every extended Young diagram we associate a \emph{charge}, $i\in I$.  In each box, we put a \emph{color} $j\in I$ given by $j \equiv a-b+i \pmod{n}$ where $a$ is the number of columns from the right and $b$ is the number of rows from the top (see figure \ref{pattern}).  

\begin{figure}[bt]
%\begin{figure}
%\caption{Color pattern for an extended Young diagram of charge $i$.  All labels are reduced modulo $n$.} \label{pattern}
\centering
\begin{tabular}{|c|c|c|c}
\hline
$i$ & $i+1$ & $i+2$ & $\cdots$\\
\hline
$i-1$ & $i$ & $i+1$ & $\cdots$\\
\hline
$i-2$ & $i-1$ & $i$ &  $\cdots$\\
\hline
\vdots & \vdots & \vdots \\
\end{tabular}
\caption{Color pattern for an extended Young diagram of charge $i$.  All labels are reduced modulo $n$.} \label{pattern}
\end{figure}

For example, {\scriptsize\young(120,01)} is an extended Young diagram of charge $1$ for $n=3$. The \emph{null diagram} with no boxes---denoted by $\varnothing$---is also considered as an extended Young diagram.

A column in an extended Young diagram is $i$-\emph{removable} if the bottom box contains $i$ and can be removed leaving another extended Young diagram.  A column is $i$-\emph{admissible} if a box containing $i$ could be added to give another extended Young diagram.

An extended Young diagram is called $n$-\emph{regular} if there are at most $(n-1)$ rows with the same number of boxes.   Let $\mathcal{Y}(0)$ denote the collection of all $n$-regular extended Young diagrams of charge $0$.  Then $\mathcal{Y}(0)$ can be given the structure of a crystal with the following actions of $\tilde{e}_i$, $\tilde{f}_i$, $\varepsilon_i,$ $\varphi_i$, and wt$(\cdot)$.  For each $i\in I$ and $b\in \mathcal{Y}(0)$ we define the \emph{$i$-signature} of $b$ to be the string of $+$'s, and $-$'s in which each $i$-admissible column receives a $+$ and each $i$-removable column receives a $-$ reading from right to left.  The \emph{reduced $i$-signature} is the result of recursively canceling all `$+-$' pairs in the $i$-signature leaving a string of the form $(-,\dots, -,+\dots, +)$.  The Kashiwara operator $\tilde{e}_i$ acts on $b$ by removing the box corresponding to the rightmost $-$, or maps $b$ to $0$ if there are no minus signs.  Similarly, $\tilde{f}_i$ adds a box to the bottom of the column corresponding to the leftmost $+$, or maps $b$ to $0$ if there are no plus signs.  The function $\varphi_i(b)$ is the number of $+$ signs in the reduced $i$-signature of $b$ and $\varepsilon_i(b)$ is the number of $-$ signs.  We define 
wt$:\mathcal{Y}(0)\rightarrow P$ by $b \mapsto \Lambda_0-\sum_{i=0}^{n-1}$\# $\{i$-colored boxes in $b\}\alpha_i$.  Then $\mathcal{Y}(0)\cong \mathcal{B}(\Lambda_0).$

\section{Decomposition of $V(\Lambda_0)^{\otimes 2}$}

In order to obtain the decomposition of $V(\Lambda_0)\otimes V(\Lambda_0)$, it suffices to find the set of the \emph{maximal} elements of the crystal base $\mathcal{B}(\Lambda_0) \otimes \mathcal{B}(\Lambda_0)$, i.e. the set of all $b_1\otimes b_2\in \mathcal{B}(\Lambda_0) \otimes \mathcal{B}(\Lambda_0)$ for which $\tilde{e}_i (b_1\otimes b_2)=0$ for all $i\in I$.  Maximal elements are characterized by the following:
\begin{lemma}{\cite[Lemma 5.1]{JMMO}}
An element $b_1 \otimes b_2 \in \mathcal{B}(\Lambda_0)\otimes \mathcal{B}(\Lambda_0)$ is maximal if and only if 
$\tilde{e}_ib_1=0 $ and $\tilde{e}_i^{\delta_{i0}+1}b_2=0$ for all $i\in I$.  
%Therefore, in particular, $b_1$ is maximal in $\mathcal{B}(\Lambda_0)$.
\end{lemma}

\begin{figure}[tb]
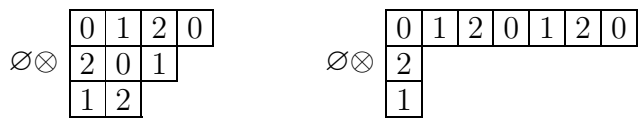

$\varnothing \otimes$ \parbox[c]{1in}{\young(0120,201,12)}\qquad $\varnothing \otimes$ \parbox[c]{2in}{\young(0120120,2,1)}
\caption{Examples of maximal elements in $\mathcal{B}(\Lambda_0)\otimes\mathcal{B}(\Lambda_0)$, for $n=3$.}  \label{delta2}
\end{figure}

We now prove the following:
\begin{lemma}\label{l2}
An element $b_1 \otimes b_2$ of the $U_q(\widehat{\mathfrak{sl}}(n))$ crystal $\mathcal{B}(\Lambda_0) \otimes \mathcal{B}(\Lambda_0)$ is maximal if and only if $b_1$ is the null diagram and the following two conditions are satisfied for $b_2$:
\begin{enumerate}
\item the first removable column from the right in $b_2$ is $0$-removable,
\item for all relevant $i$ if the $i$th admissible column in $b_2$ is $j$-admissible then the $i+1$st removable column, if it exists,  is $j$-removable for some $j \in \mathbb{Z}/n\mathbb{Z}$.
\end{enumerate}
\end{lemma} 
\begin{proof}
	It is easy to see that conditions (1) and (2) are sufficient to conclude that $b_1 \otimes b_2$ is maximal.  To show the necessity, suppose $b_1 \otimes b_2$ is a maximal element of $\mathcal{B}(\Lambda_0) \otimes \mathcal{B}(\Lambda_0)$.  To prove condition (1), let the first removable column from the right in $b_2$ be $j$-removable.  Therefore, there is a $-$ in the $j$-signature of $b_2$ and no $+$ to its left, since $b_2$ is $n$-regular.  Therefore $\varepsilon_j(b_2)\geq 1$, so $j=0$ by lemma 3.1.  To prove condition (2), let $i$ be the least index such that the $i$th admissible column from the right in $b_2$ is $j$-removable and the $i+1$st removable column is $k$-removable and $k\neq j$.  There are two cases to consider: $k\neq 0$ or $k=0$.  

\emph{Case 1:} Suppose that $k\neq 0$.  Consider the $k$-signature of $b_2$.  If there is a $+$ coming from some admissible column, say the $l$th where $l<i$, then by the minimality of $i$ it must be followed by a $-$ coming from the $l+1$st admissible column.  However, there can be nothing in between. If there were then it must come from the $l$th removable column which cannot also be $k$-removable since $b_2$ is $n$-regular.  Therefore, the $k$-signature alternates between $+$ and $-$ up to the $i+1$st removable column.  Therefore, in the reduced $k$-signature the $-$ coming from the $i+1$st removable column cannot be canceled, contradicting the maximality of $b$.  Therefore, there is no such index $i$ and we are done.

\emph{Case 2:} Suppose that $k=0$.  Since $\varepsilon_0(b_2)=1$ there must be exactly one $-$ in the reduced 0-signature of $b_2$.  However, the first removable column already contributes one $-$ sign, so the same argument as case 1 shows that we can cancel out every $-$ after the first $-$ up to the $i+1$st removable column.  Therefore no such index $i$ exists.
\end{proof}
As an application of lemma \ref{l2}, all maximal elements of weight $2\Lambda_0-3\delta$ for $\mathcal{B}(\Lambda_0)\otimes\mathcal{B}(\Lambda_0)$ with $n=3$ are given in figure \ref{delta2}.\\
As usual, we denote a \emph{partition} by a finite sequence $(\lambda_1^{f_1},\lambda_2^{f_2},\dots, \lambda_j^{f_j}),$ where $\lambda_k\in \mathbb{Z}_{> 0},$ $\lambda_k>\lambda_{k+1},$ and $f_k\in \mathbb{Z}_{> 0}$ denotes the multiplicity of 
$\lambda_k$.  Each $b\in\mathcal{\mathcal{Y}}(0)$ can be uniquely represented as a partition where $\lambda_k$ is the number of boxes in a given row, and $f_k$ is the number of rows having $\lambda_k$ boxes.  For example, the two diagrams in figure \ref{delta2} correspond to the partitions $(4,3,2),$ and $(7,1^2).$  We can now rephrase lemma 3.2 in terms of partitions as follows.
\begin{lemma}
The maximal elements of $\mathcal{B}(\Lambda_0)\otimes\mathcal{B}(\Lambda_0)$ are in a one-to-one correspondence with the set of all partitions $\lambda=(\lambda_1^{f_1},\lambda_2^{f_2},\dots, \lambda_j^{f_j})$ such that $f_k< n,$ $k=1,2,\dots, j$,  and satisfying the conditions:

\begin{enumerate}
\item $f_1 \equiv \lambda_1 \pmod{n}$,
\item $f_k+f_{k+1}+\lambda_{k}-\lambda_{k+1}\equiv0 \pmod{n}$, for $k<j$.
\end{enumerate}
\end{lemma}

\begin{proof} The condition that each $f_k< n$ is equivalent to the condition that $b\in\mathcal{\mathcal{Y}}(0)$ is $n$-regular.  
The first column from the right is always removable, so by condition (1) of lemma 3.2 the first column contains a 0-colored box  in the bottom row.  Since $\lambda_1$ is the number of columns in the diagram and $f_1$ is the number of boxes in the rightmost column, we see that $\lambda_1-f_1\equiv 0 \pmod{n}$.

Now, suppose that the $k$th admissible column is $t$-admissible, and there exists a removable column to the left of that column.  The $k$th admissible column is $\lambda_{k}+1$ columns from the left and contains $f_1+f_2+\dots +f_{k-1}$ boxes.  Therefore: 
\begin{eqnarray}
f_1+f_2+\dots +f_{k-1}-(\lambda_k+1)&\equiv&t-1 \pmod{n}\nonumber\\
f_1+f_2+\dots +f_{k-1}-\lambda_k&\equiv&t \pmod{n}\label{t1}
\end{eqnarray}

By condition (2) of lemma 3.2 the $k+1$st removable column (the $\lambda_{k+1}$st from the left) is also $t$-removable, and contains $f_1+f_2+\dots +f_k+f_{k+1}$ boxes.   Therefore:
\begin{equation}
f_1+f_2+\dots+f_k+f_{k+1}-\lambda_{k+1} \equiv t \pmod{n}\label{t2}
\end{equation}
Subtracting equation \eqref{t1} from equation \eqref{t2} we obtain condition (2). Furthermore, if the partition satisfies the conditions  (1) and (2), then it is in correspondence with an extended Young diagram as in Lemma 3.2.
\end{proof}
Let $\mathscr{C}_n$ be the collection of all partitions satisfying conditions (1) and (2) in lemma 3.3. Then each such partition in $\mathscr{C}_n$ corresponds to a unique maximal element in  $\mathcal{B}(\Lambda_0)\otimes\mathcal{B}(\Lambda_0)$. Let  
$\mathscr{C}_{n,\mu}$ denote the set of elements in $\mathscr{C}_n$ of weight $\mu$. It is known that each connected component of the crystal $\mathcal{B}(\Lambda_0)\otimes\mathcal{B}(\Lambda_0)$ is  the crystal for the corresponding irreducible summand of the module $V(\Lambda_0)\otimes V(\Lambda_0)$. In the following lemma we determine the connected components of the crystal $\mathcal{B}(\Lambda_0)\otimes\mathcal{B}(\Lambda_0)$.

\begin{lemma}
Each connected component of $\mathcal{B}(\Lambda_0)\otimes\mathcal{B}(\Lambda_0)$ is isomorphic to \\
$\mathcal{B}\left (\Lambda_i+\Lambda_{n-i}-k\delta\right )$ for some $i\in\{0,1,\dots, [n/2]\}$ and $k\in \mathbb{Z}_{\geq 0}$ such that $k\geq i$.
\end{lemma}
\begin{proof} 
Let $b\otimes b'\in \mathcal{B}(\Lambda_0)\otimes\mathcal{B}(\Lambda_0)$ be maximal.  Then $\text{wt}(b)=\Lambda_0$, and $b'$ corresponds to a partition $(\lambda_1^{f_1},\lambda_2^{f_2}\dots,\lambda_l^{f_l})\in \mathscr{C}_n$.  We set $\alpha=\Lambda_0-\text{wt}(b)$ and use the weight formula in $\mathcal{B}(\Lambda_0)$  to compute $\alpha$:

\begin{eqnarray}
\alpha&=&\sum_{t=1}^{l}\sum_{u=1}^{\lambda_{t}}\sum_{v=s_{t-1}+1}^{s_{t}}\alpha_{u-v}, \label{trickysum}\\
	&=&\sum_{t=1}^{l}\sum_{u=1}^{\lambda_{t}}\sum_{v=s_{t-1}+1}^{s_{t}}(2\Lambda_{u-v}-\Lambda_{u-v-1}-\Lambda_{u-v+1}+\delta_{\overline{u-v},0}\delta),
\end{eqnarray}
where $ s_t:=\sum_{m=1}^t f_m$.  The sums telescope, leaving:

\begin{eqnarray}
\alpha&=&
{\sum_{t=1}^l(\Lambda_{\lambda_t-s_t}-\Lambda_{-s_t}+\Lambda_{-s_{t-1}}-
\Lambda_{\lambda_t-s_{t-1}})}+k\delta,\label{withdelta}
\end{eqnarray}
where $k\in \mathbb{Z}_{\geq 0}$ is defined to be the number of $\alpha_0$'s in the sum \eqref{trickysum}.  Now consider the sum in \eqref{withdelta}.

\qquad\begin{eqnarray*}
\lefteqn{\sum_{t=1}^l(\Lambda_{\lambda_t-s_t}-\Lambda_{-s_t}+\Lambda_{-s_{t-1}}-
\Lambda_{\lambda_t-s_{t-1}})}\\
	&=&\Lambda_{\lambda_1-f_1}+\sum_{t=2}^l\Lambda_{\lambda_{t-1}+f_{t-1}+f_t-s_{t}}-\sum_{t=1}^l\Lambda_{-s_{t}}
	+\sum_{t=1}^l\Lambda_{-s_{t-1}}-\sum_{t=1}^l\Lambda_{\lambda_t-s_{t-1}}\\
	&&\qquad \text{since }\lambda_t\equiv \lambda_{t-1}+f_{t-1}+f_{t}\pmod{n} \\
&=&\Lambda_{0}+\sum_{t=2}^l\Lambda_{\lambda_{t-1}-s_{t-2}}-\sum_{t=1}^l
	\Lambda_{\lambda_t-s_{t-1}}+\sum_{t=1}^l\Lambda_{-s_{t-1}}
	-\sum_{t=1}^l\Lambda_{-s_{t}}\\
	&&\qquad \text{since }\lambda_1\equiv f_1 \pmod{n}\\
&=&\Lambda_{0}+\sum_{t=1}^{l-1}\Lambda_{\lambda_{t}-s_{t-1}}
	-\sum_{t=1}^l\Lambda_{\lambda_t-s_{t-1}}
	+\sum_{t=0}^{l-1}\Lambda_{-s_{t}}-\sum_{t=1}^l\Lambda_{-s_{t}}\\
&=&\Lambda_0-\Lambda_{\lambda_l-s_{l-1}}+\Lambda_0-\Lambda_{-s_l}.
\end{eqnarray*}

Let $i:=\min\{\lambda_l-\overline{s_{l-1}}, \overline{-s_l}\}$ and $j:=\max\{\lambda_l-\overline{s_{l-1}},\overline{-s_l}\}$.  By repeated use of condition (2) of Lemma 3.3 we can see that $i+j\equiv \lambda_l-s_{l-1}-s_l\equiv \lambda_1-s_0-s_1=\lambda_1-f_1 \equiv0 \pmod{n}.$  Since  $i\le j$ we have that $i\in\{0,1,\dots,\left [n/2\right ]\}$ and $j=n-i$.  Here and after $[\cdot]$ denotes the floor function.  

Since $\Lambda_t=\Lambda_0+\omega_t$ (\cite{K} , eq. 12.4.3), using  well-known formulas for the fundamental dominant weights $\omega_t$ of $\mathfrak{sl}(n)$ (see for example \cite{H}) we have:
\begin {eqnarray*}
\alpha&=&2\Lambda_0-\Lambda_i-\Lambda_{n-i}+k\delta\\
	&=&2\Lambda_0-\Lambda_0-\frac{1}{n}\left (\sum_{t=1}^{i}t(n-i)\alpha_t+
		\sum_{t=i+1}^{n-1}i(n-t)\alpha_t\right )\\
	&&-\:\Lambda_0-\frac{1}{n}\left (\sum_{t=1}^{n-i}ti\alpha_t+
		\sum_{t=n+1-i}^{n-1}(n-i)(n-t)\alpha_t\right )+k\delta\\
	&=&-\sum_{t=1}^{i}t\alpha_t-\sum_{t=i+1}^{n-i}
	i\alpha_t-\sum_{t=n+1-i}^{n-1}(n-t)\alpha_t+k\delta\\
	&=&i\alpha_0+\sum_{t=1}^{i-1}(i-t)(\alpha_t+\alpha_{n-t})+(k-i)\delta.
\end{eqnarray*}
Because $\alpha\in Q^+$ we must have  $k\geq i$, otherwise the coefficient of $\alpha_{[n/2]}$ would be $<0$ since $i\leq [n/2]$.  Therefore $\text{wt}(b\otimes b')=2\Lambda_0-\alpha=\Lambda_i+\Lambda_{n-i}-k\delta.$  Hence $b\otimes b'$ is a maximal element of weight $\Lambda_i+\Lambda_{n-i}-k\delta$, hence the component of $\mathcal{B}(\Lambda_0)\otimes\mathcal{B}(\Lambda_0)$ containing $b\otimes b'$ is isomorphic to $\mathcal{B}(\Lambda_i+\Lambda_{n-i}-k\delta)$ with $k\geq i.$
\end{proof}

Now, consider the case when $k=i$.  Note that the coefficient of $\alpha_0$ is $i$.  If two `0' s appeared in the same row or column, then the coefficient of $\alpha_t$ in $\alpha$ would be $>0$ for all $t\in I$, which cannot happen.  Also, $\text{ht}(\alpha)=i+\sum_{t=1}^{i-1}2(i-t)=i^2,$ so the corresponding partition must have $i^2$ boxes.  Therefore, the only maximal element of weight $\Lambda_i+\Lambda_{n-i}-i\delta$ is a square partition with side $i$ (see figure \ref{squares}).

\begin{figure}
%\caption{Maximal elements of weight $2\Lambda_0-\alpha_{tn}$ for $n=6$ and $t=0,1,2,3$.}
%\label{squares}
$\varnothing \otimes\varnothing$ \qquad $\varnothing \otimes$ \parbox[c]{1in}{\young(0)}$\varnothing \otimes$ \parbox[c]{1in}{\young(01,50)}\qquad$\varnothing \otimes$ \parbox[c]{1in}{\young(012,501,450)}
\caption{Maximal elements of weight $\Lambda_i+\Lambda_{n-i}-i\delta$ for $n=6$ and $i=0,1,2,3$.}
\label{squares}
\end{figure}

Using the above Lemma we have now proved our main result:

\begin{theorem}
The $\widehat{\mathfrak{sl}}(n)$-module $V(\Lambda_0)\otimes V(\Lambda_0)$ decomposes as the direct sum $\bigoplus_{i=0}^{[n/2]} \bigoplus_{k=i}^{\infty} V(\Lambda_i+\Lambda_{n-i}-k\delta)^{\oplus b_{ik}}$ and the outer multiplicities are given by $$b_{ik}=|\mathscr{C}_{n,\Lambda_i+\Lambda_{n-i}-k\delta}|,$$
where the absolute value sign denotes the cardinality.
Here, if $b_{ik}=0$ then $V(\Lambda_i+\Lambda_{n-i}-k\delta)$ does not occur in the decomposition.  Moreover, $b_{ii}=1$.
\end{theorem}

As an example of Theorem 3.1 consider the decomposition of the $\widehat{\mathfrak{sl}}(3)$-module $V(\Lambda_0)^{\otimes 2}\cong \bigoplus_{k=0}^{\infty}\left (V(2\Lambda_0-k\delta)^{\oplus b_{0k}}\oplus V(2\Lambda_0-\alpha_0-k\delta)^{\oplus b_{1k}}\right )$. The following table gives the outer multiplicities for $k=0,\dots, 6:$

\begin{center}
\begin{tabular}{|c|c|c|c|c|c|}
\hline
$k$ & $\mathscr{C}_{3,2\Lambda_0-k\delta}$ & $b_{0,k}$& $k$ & $\mathscr{C}_{3,2\Lambda_1-k\delta}$ &  $b_{1,k}$\\
\hline
0 & $()$ & 1 &1& $(1)$ & 1\\
\hline
1 & & 0&2 & $(4), (2^2)$ & 2\\
\hline
2 & $(4,1^2)$ & 1&3 & $(7), (4,3)$ & 2\\
\hline
3& $(7,1^2),(4,3,2)$ & 2&4 & $(10),(7,3),(5^2),(4,3,2,1)$ & 4\\
\hline
4 & $(10,1^2), (7,3,2), (5^2,2)$ & 3&5 & $(13), (10,3), (7,6)$ &5 \\& & && $(7,3,2,1), (5^2,2,1)$ & \\
\hline
5 & $(13,1^2), (10,3,2), (7,6,2), (7,4^2)$ & 4&6 &$ (16), (13,3), (10,6) $& 8\\

&  & && $ (10,3,2,1),(8^2), (7,6,2,1)$ &\\
&&& &$(7,4^2, 1), (5^2, 3^2)$ &\\
\hline
6 & $(16, 1^2), (13,3,2), (10,6,2), (10,4^2)$ & 7&7 & $(19), (16,3), (13,6)$ & 11\\
& $(8^2,2), (7,6,5), (5^2, 3^2, 1^2)$ & && $ (13,3,2,1),(10,9), (10,6,2,1)$ & \\
& & & & $ (10,4^2,1),(8^2,2,1),(7,6,5,1)$ &\\
& & & &$(7,6,3^2),(7,4^2, 2^2)$ &\\ 
\hline
\end{tabular}
\end{center}
\section{Generating Functions for Outer Multiplicities}
In this section we consider the generating functions for the outer multiplicities $b_{ik}$ in Theorem 3.1 and give explicit formulas for these generating functions.  First we define  the formal series $f$ and $g$  in the indeterminates $u,v$ as follows: 
\begin{eqnarray*}
f(u,v)&=&\sum_{j=-\infty}^\infty u^{j(j-1)/2}v^{j(j+1)/2},\\
g(u,v)&=&\sum_{j=-\infty}^\infty (-1)^j u^{j(j-1)/2}v^{j(j+1)/2}.
\end{eqnarray*}
It is easy to verify that the functions $f(u,v)$ and $g(u,v)$ satisfy the following properties: 

\begin{eqnarray}
f(u,v) = f (v,u), \quad \text{and} \quad g(u,v) = g(v,u),\\
f(q^r,q^s) = q^r f(q^{2r+s},q^{-r}), \quad \text{and} \quad g(q^r,q^s)&=&-q^r g(q^{2r+s},q^{-r}),
\end{eqnarray}

\noindent for integers $r, s$ not both equal to $0$.
Recall the Euler $\varphi$ function  $\varphi(q):=g(q,q^2)=\prod_{j=1}^\infty (1-q^j)$. 
In what follows we will be using the well known Jacobi triple product identities (cf. \cite{A}):
\begin{eqnarray}
f(u,v)&=&\prod_{j=1}^\infty (1-u^j v^j)(1+u^{j-1} v^j)(1+u^j v^{j-1}),\\
g(u,v)&=&\prod_{j=1}^\infty (1-u^j v^j)(1-u^{j-1} v^j)(1-u^j v^{j-1}).
\end{eqnarray}

Recall the $q$-character formulas (or principally specialized character formulas)  for the  $\widehat{\mathfrak{sl}}(n)$-modules $V(\Lambda_i), i \in \{0, 1, \dots , n-1\}$ and $V(\Lambda_0+\Lambda_j)$ for $j\in \{ 0,1, \dots, \left [ n/2 \right ]\}$ (cf. \cite{M}):

\begin{equation}
\text{ch}_q(V(\Lambda_i))=\prod_{\substack{j>0 \\ j\not \equiv 0 \pmod{n}}}(1-q^j)^{-1}= \frac{\varphi(q^{n})}{\varphi(q)}.
\label{qdim1}
\end{equation}
\begin{equation}
\text{ch}_q(V(\Lambda_0+\Lambda_j))=\frac{\varphi(q^{n})g(q^{j+1},q^{n-j+1})}{\varphi(q)^2}.
\label{qdim2}
\end{equation}
  
Now we define the generating function
\begin{equation}
B_i(q)=\sum_{k=i}^\infty b_{ik} q^{k-i} ,
\end{equation}
for the outer multiplicities $b_{ik}$.  Note that $2\Lambda_0=\Lambda_i+\Lambda_{n-i}+i\alpha_0+\sum_{t=1}^{i-1}(i-t)(\alpha_t+\alpha_{n-t})-i\delta$.  By Theorem 3.1 we have:
\begin{eqnarray*}
\text{ch}(V(\Lambda_0)\otimes V(\Lambda_0))&=& \sum_{i=0}^{[n/2]} \sum_{k=i}^{\infty} b_{ik}\:\text{ch}(V(\Lambda_i+\Lambda_{n-i}-k\delta).
\end{eqnarray*}

Hence 
\begin{eqnarray*}
e(-2\Lambda_0)\text{ch}(V(\Lambda_0))^2&=& e(-2\Lambda_0)\sum_{i=0}^{[n/2]} \sum_{k=i}^{\infty} b_{ik}\: \text{ch}(V(\Lambda_i+\Lambda_{n-i}))e(-k\delta)\\
	&=& \sum_{i=0}^{[n/2]} \sum_{k=i}^{\infty} \bigg [b_{ik}\: e(-\Lambda_i-\Lambda_{n-i})\\
	&&\cdot\:e\left (-i\alpha_0-\sum_{t=1}^{i-1}(i-t)(\alpha_t+\alpha_{n-t})+i\delta\right )\\
	&&\cdot \: \text{ch}(V(\Lambda_{i}+\Lambda_{n-i})e(-k\delta) \bigg ]\\
	&=& \sum_{i=0}^{\left [ n/2\right ]}\bigg [e\left (-i\alpha_0-\sum_{t=1}^{i-1}(i-t)
	(\alpha_t+\alpha_{n-t})\right )\\
	&& \cdot \: e(-\Lambda_i-\Lambda_{n-i})\text{ch}(V(\Lambda_i+\Lambda_{n-i}))\sum_{k=i}^\infty  b_{ik}e(-(k-i)\delta)\bigg ].
\end{eqnarray*}
Now specializing $e(-\alpha_i)=q$, we obtain:

\begin{eqnarray*}
\text{ch}_q(V(\Lambda_0))^2&=&\sum_{i=0}^{\left [ n/2 \right ]} q^{i^2}\text{ch}_q(V(\Lambda_i+\Lambda_{n-i}))\sum_{k=i}^\infty b_{ik} q^{n(k-i)},
\end{eqnarray*}

which gives

\begin{eqnarray*}
\frac{\varphi(q^{n})^2}{\varphi(q)^2}&=&\sum_{i=0}^{\left [ n/2\right ]}q^{i^2}\frac{\varphi(q^{n})g(q^{2i+1},q^{n+1-2i})}{\varphi(q)^2}B_i(q^{n}),
\end{eqnarray*}
and hence we have
\begin{equation}
\varphi(q^{n})=\sum_{i=0}^{\left [ n/2 \right ]}q^{i^2}g(q^{2i+1},q^{n+1-2i})B_i(q^{n})\label{eq2}.
\end{equation}

Note that the series $\varphi(q^{n})=\prod_{j=1}^\infty(1-q^{nj})$ has a zero coefficient in front of $q^j$ whenever $j$ is not a multiple of $n$, and similar is the case for  $B_i(q^{n})$.  However, this is not the case for $q^{i^2}g(q^{2i+1},q^{n+1-2i})$.   So the trick is to rearrange the sum to sort the powers of $q$ carefully as we do  below.\\

Set  $\Phi_{in}(q):=q^{i^2}g(q^{2i+1},q^{n+1-2i})$.  Then:
\begin{eqnarray*}
\Phi_{in}(q)&=&\sum_{j=0}^{n-1} \sum_{\substack{m\in \mathbb{Z}\\ m\equiv j \pmod{n}}}(-1)^j q^{j[(n+2)j+4i-n]/2+i^2}\\
	&=&\sum_{j=0}^{n-1}\sum_{m \in \mathbb{Z}}(-1)^{nm+j}q^{(nm+j)[(n+2)(nm+j)+4i-n]/2+i^2}\\
	&=&\sum_{j=0}^{n-1} (-1)^j q^{ nj(j-1)/2+(i+j)^2}h\left (q^{n[n(n+3)/2-2i-(n+2)j ]},q^{n[n(n+1)/2+2i+(n+2)j)]}\right )
\end{eqnarray*}
where $h=f$ if $n$ is even and $h=g$ if $n$ is odd.  Define
\begin{equation*}
\Psi_{ijn}(q):=h\left (q^{n(n+3)/2-2i-(n+2)j},q^{n(n+1)/2+2i+(n+2)j}\right ).
\end{equation*}
If $n$ is even \eqref{eq2} becomes:
\begin{eqnarray}
\varphi(q^{n})&=&\sum_{i=0}^{n/2}\sum_{j=0}^{n-1} (-1)^j q^{ nj(j-1)/2+(i+j)^2}\Psi_{ijn}(q^{n})B_i(q^{n})\nonumber\\
	&=&\sum_{i=0}^{n/2}B_i(q^{n})\big [ (-1)^{\overline{-i}} q^{ n[\overline{-i}(\overline{-i}-1)]/2+(\overline{-i}+i)^2}\Psi_{i,\overline{-i},n}(q^{n})\nonumber\\
	&&+\:(-1)^{n/2-i} q^{ n(n/2-i)(n/2-i-1)/2+(n/2)^2}\Psi_{i,n/2-i,n}(q^{n})\nonumber\\
	&&+\:\sum_{j=1}^{(n-2)/2}\big\{(-1)^{j-i}q^{ n[\overline{j-i}(\overline{j-i}-1)]/2+(\overline{j-i}+i)^2}\Psi_{i,\overline{j-i},n}(q^{n})\nonumber\\
	&&+\:(-1)^{i+j} q^{n[\overline{-i-j}(\overline{-i-j}-1)]/2+(\overline{-i-j}+i)^2}\Psi_{i,\overline{-i-j},n}(q^{n})\big \} \big ],\label{lodd}
\end{eqnarray}
If $n$ is odd, then \eqref{eq2} becomes: 
\begin{eqnarray}
\varphi(q^{n})&=&\sum_{i=0}^{(n-1)/2}\sum_{j=0}^{n-1} (-1)^j q^{ n j(j-1)/2+(i+j)^2}\Psi_{ijn}(q^{n})B_i(q^{n})\nonumber\\
	&=&\sum_{i=0}^{(n-1)/2}B_i(q^{n})\big [ (-1)^{\overline{-i}} q^{ n[\overline{-i}(\overline{-i}-1)]/2+(\overline{-i}+i)^2}\Psi_{i,\overline{-i},n}(q^{n})\nonumber\\
	&&+\:\sum_{j=1}^{(n-1)/2}\big\{(-1)^{\overline{j-i}}q^{ n[\overline{j-i}(\overline{j-i}-1)]/2+(\overline{j-i}+i)^2}\Psi_{i,\overline{j-i},n}(q^{n})\nonumber\\
	&&+\:(-1)^{\overline{-i-j}} q^{n[\overline{-i-j}(\overline{-i-j}-1)]/2+(\overline{-i-j}+i)^2}\Psi_{i,\overline{-i-j},n}(q^{n})\big \} \big ].\label{leven}
\end{eqnarray}

Notice that in each term of the inner summations of \eqref{lodd} and \eqref{leven}, each of the exponents of $q$ in both summands has the same remainder modulo $n$---namely $j^2 \mod{n}$---while the first term in the outer summation contains exponents of $q$ that reduce to $0 \mod{n}.$  Now, suppose that  $n=p$ is an odd prime so that all the quadratic residues mod $n$ are distinct.  Then we are justified in separating \eqref{leven} into the following $(p+1)/2$ equations with the $(p+1)/2$ generating functions as unknowns.
\begin{eqnarray*}
\sum_{i=0}^{(p-1)/2}  (-1)^{\overline{-i}} q^{ p[\overline{-i}(\overline{-i}+1)]/2+(\overline{-i}+i)^2}\Psi_{i,\overline{-i},p}(q^{p})B_i(q^{p})&=&\varphi(q^{p}),\\
	\sum_{i=0}^{(p-1)/2}\big \{(-1)^{\overline{j-i}}q^{ p[\overline{j-i}(\overline{j-i}-1)]/2+(\overline{j-i}+i)^2}\Psi_{i,\overline{j-
	i},p}	(q^{p})\qquad \qquad\quad\\
	+\:(-1)^{\overline{-i-j}} q^{p[\overline{-i-j}(\overline{-i-j}-1)]/2+(\overline{-i-j}+i)^2}\Psi_{i,\overline{-i-j},p}(q^{p})
	\big \}B_i(q^{p})&=&0,
\end{eqnarray*}
where $j=1,\dots,(p-1)/2$.  We can factor and cancel the remainders modulo $p$.  This leaves only powers of $q^{p},$ so we may let $q\mapsto q^{1/p}$ to get:
\begin{eqnarray}
\sum_{i=0}^{(p-1)/2}  (-1)^{\overline{-i}} q^{\overline{-i}(\overline{-i}-1)/2+(\overline{-i}+i)^2}\Psi_{i,\overline{-i},p}(q)B_i(q)&=&\varphi(q),\label{gfeqn1}\\
	\sum_{i=0}^{(p-1)/2}\big \{(-1)^{\overline{j-i}}q^{ \overline{j-i}(\overline{j-i})/2+Q([\overline{j-i}+i]^2)}\Psi_{i,\overline{j-
	i},p}(q)\qquad \qquad\quad\label{gfeqn2}\\ 
	+\:(-1)^{\overline{-i-j}} q^{\overline{-i-j}(\overline{-i-j}-1)/2+Q([\overline{-i-j}+i]^2)}\Psi_{i,\overline{-i-j},p}(q)
	\big \}B_i(q)&=&0,\nonumber
\end{eqnarray}
where $Q(t)$ denotes the integer quotient function of $t$ by $p,$ ie. the unique integer such that $t=pQ(t)+\overline{t}.$  Now we can solve the following matrix equation:
\begin{equation}
A\begin{pmatrix}
		B_0(q)\\
		B_1(q)\\
		\vdots\\
		B_{(p-1)/2}(q)
	      \end{pmatrix}	
		=\begin{pmatrix}\varphi(q)\\
				0\\
				\vdots\\
				0
		\end{pmatrix},\label{mateqn}	
\end{equation}
where $A=\left ( a_{ij}\right )_{(p+1)/2 \times (p+1)/2}$ is the coefficient matrix of \eqref{gfeqn1} and \eqref{gfeqn2}, with both indices ranging from $0$ to $(p-1)/2$.  Explicitly:
\begin{eqnarray*}
a_{0i}&=&(-1)^{\overline{-i}} q^{\overline{-i}(\overline{-i}-1)/2+Q(\overline{-i}+i)^2}\Psi_{i,\overline{-i},p}(q),\\
a_{ji}&=&(-1)^{\overline{j-i}}q^{ \overline{j-i}(\overline{j-i}-1)/2+Q([\overline{j-i}+i]^2)}\Psi_{i,\overline{j-
	i},p}(q)\qquad \qquad\quad\\ 
	&&+\:(-1)^{\overline{-i-j}} q^{\overline{-i-j}(\overline{-i-j}-1)/2+Q([\overline{-i-j}+i]^2)}\Psi_{i,\overline{-i-j},p}(q),\qquad j\neq 0.
\end{eqnarray*}

Then, using Cramer's rule in \eqref{mateqn} gives the following expression for $B_i(q):$
\begin{equation}
B_i(q)=\frac{(-1)^i\varphi(q)\det(\widetilde{A_{0i}})}{\det(A)},\label{deteqn}
\end{equation}
where $\widetilde{A_{ij}}$ denotes the matrix $A$ with row $i$ and column $j$ deleted.  Recall that we assume here that $p$ is an odd prime, but we conjecture that \eqref{deteqn} holds for all odd $n$.

In a similar computation to the case where $n$ is an odd prime, we can show that \eqref{deteqn}  holds for $n=2p$ where $p=1$ or  an odd prime, where $A=\left ( a_{ij}\right )_{p+1 \times p+1}$ with indices ranging from $0$ to $p$ and 

\begin{eqnarray*}
a_{0i}&=&(-1)^{\overline{-i}} q^{\overline{-i}(\overline{-i}-1)/2+Q(\overline{-i}+i)^2}\Psi_{i,\overline{-i},2p}(q),\\
a_{ji}&=&(-1)^{\overline{j-i}}q^{ \overline{j-i}(\overline{j-i}-1)/2+Q([\overline{j-i}+i]^2)}\Psi_{i,\overline{j-
	i},2p}(q)\qquad \qquad\quad\\ 
	&&+\:(-1)^{\overline{-i-j}} q^{\overline{-i-j}(\overline{-i-j}-1)/2+Q([\overline{-i-j}+i]^2)}\Psi_{i,\overline{-i-j},2p}(q),\qquad j\neq 0,p,\\
a_{pi}&=&(-1)^{p-i}q^{ (p-i)(p-i-1)/2+Q(p^2)}\Psi_{i,p-i,2p}.
\end{eqnarray*}
We also conjecture that \eqref{deteqn} holds for all even $n$ with $A$ as given above.

\section{Some Examples and Generating Function Identities}
In this section we consider the cases $n = 2$, and $3$. Using the two different ways of computing the outer multiplicities, we obtain some interesting  generating function identities.  For $n=2,$ we have $\mathscr{C}_2=\{\text{all partitions with distinct odd parts}\}$, $\mathscr{C}_{2,2\Lambda_0-k\delta}=\{\text{all partitions of }2k$ with distinct odd parts$\}$,
$\mathscr{C}_{2,2\Lambda_1-k\delta}=\{\text{all partitions of }2k-1$ with distinct odd parts$\}$.  Therefore, by Theorem 3.1 using elementary generating function methods (see \cite{A} for example) we obtain following  infinite sum forms for the generating functions of the outer multiplicities: 
\begin{eqnarray}
B_1(q)&=&\sum_{i=0}^{\infty}\frac{q^{2i^2}}{\prod_{k=1}^{2i}(1-q^k)},\\
B_2(q)&=&\sum_{i=0}^{\infty}\frac{q^{2i^2+2i}}{\prod_{k=1}^{2i+1}(1-q^k)}.
\end{eqnarray}

On the other hand in the case $n = 2$ it follows from  \eqref{deteqn} that 

\begin{eqnarray}
B_1(q)&=&\varphi(q)\frac{f(q^{5},q^{3})}{D_2},\\
B_2(q)&=&\varphi(q)\frac{f(q,q^{7})}{D_2},
\end{eqnarray}
where $D_2:=f(q^5,q^3)^2-qf(q,q^7)^2$.  Thus we obtain the identities:

\begin{lemma}
\begin{eqnarray}
\sum_{i=0}^{\infty}\frac{q^{2i^2}}{\prod_{k=1}^{2i}(1-q^k)}&=& \varphi(q)\frac{f(q^{5},q^{3})}{D_2},\\
\sum_{i=0}^{\infty}\frac{q^{2i^2+2i}}{\prod_{k=1}^{2i+1}(1-q^k)} &=& \varphi(q)\frac{f(q,q^{7})}{D_2}.
\end{eqnarray}
\end{lemma}

Recall the following known identities (\textbf{S.38} and \textbf{S.39} in \cite{MSZ}):

\begin{eqnarray}
\frac{f(q^5,q^3)}{\varphi(q^2)}&=&\sum_{i=0}^{\infty}\frac{q^{2i^2}}{\prod_{k=1}^{2i}(1-q^k)},\\
\frac{f(q,q^7)}{\varphi(q^2)}&=&\sum_{i=0}^{\infty}\frac{q^{2i^2+2i}}{\prod_{k=1}^{2i+1}(1-q^k)}.
\end{eqnarray}

Hence it follows from Lemma 5.1 that:

\begin{lemma}
$f(q^5,q^3)^2-qf(q,q^7)^2=\varphi(q)\varphi(q^2).$
\end{lemma}

In the $n=2$ case, Feingold \cite{F} gave the following generating functions for the outer multiplicities:
\begin{eqnarray}
B_1(q)&=&\frac{f(q^{11},q^{13})-qf(q^5,q^{19})}{\varphi(q)},\\
B_2(q)&=&\frac{f(q^{7},q^{17})-q^2f(q,q^{23})}{\varphi(q)}.
\end{eqnarray}

Comparing these with our result in (5.3) and (5.4) and Lemma 5.2 we obtain the identities:

\begin{lemma}
\begin{eqnarray}
\frac{f(q^{5},q^{3})}{\varphi(q^2)}&=&\frac{f(q^{11},q^{13})-qf(q^5,q^{19})}{\varphi(q)},\\
\frac{f(q,q^{7})}{\varphi(q^2)}&=&\frac{f(q^{7},q^{17})-q^2f(q,q^{23})}{\varphi(q)}.
\end{eqnarray}
\end{lemma}

\begin{proof}
It follows from Lemma 5.1, and Lemma 5.2 using the known identities (\textbf{S.83}) and (\textbf{S.86})  in  \cite{MSZ}.
\end{proof}

Now we turn to the case $n=3$.  Since $n$ is an odd prime, \eqref{deteqn} applies and yields:
\begin{eqnarray}
B_1(q)&=&\frac{\varphi(q)(g(q^7,q^8)- q\cdot g(q^{2},q^{13}))}{D_3},\\
B_2(q)&=&\frac{\varphi(q)(g(q^{11},q^4)+q\cdot g(q,q^{14}))}{D_3},
\end{eqnarray}
where
\begin{eqnarray*}
D_3&:=&g(q^6,q^9)(g(q^7,q^8)- q\cdot g(q^{2},q^{13}))-q\cdot g(q^{12},q^3)(g(q^{11},q^4)\\&&+\:q\cdot g(q,q^{14})).
\end{eqnarray*}

\begin{lemma}
$D_3=\varphi(q)^2.$
\end{lemma}
\begin{proof}  
Consider the product:
\begin{eqnarray}
\varphi(q)^2&=&g(q,q^2)^2\nonumber\\
	&=&\sum_{n=-\infty}^\infty (-1)^n q^{n(3n+1)/2}\sum_{m=-\infty}^\infty (-1)^m q^{m(3m+1)/2}\nonumber\\
	&=&\sum_{(m,n)\in \mathbb{Z}^2} (-1)^{m+n}q^{n(3n+1)/2+m(3m+1)/2}. \label{detsum}
\end{eqnarray}
We now separate the sum in \eqref{detsum} into cosets of a particular sublattice of $\mathbb{Z}^2$.  We choose the sublattice $L:=\text{span}_{\mathbb{Z}}\{(-1,2),(-2,-1)\},$ which is equivalent to making the substitution $m=-i-2j,$ $n=2i-j.$   The index of $L$ in $\mathbb{Z}^2$ is 5, and its cosets, which form a partition of $\mathbb{Z}^2$, are: $L$, $L+(-1,0)$, $L+(-2,0)$, $L+(-1,1),$ and $L+(-2,1)$.  So the sum in \eqref{detsum} becomes:
\begin{eqnarray*}
\sum_{\{M | M\text{is a coset of }L\}}\sum_{(m,n)\in M} (-1)^{m+n}q^{n(3n+1)/2+m(3m+1)/2}.
\end{eqnarray*}
This simplifies to:
\begin{eqnarray}
\varphi(q)^2 &=&g(q^7,q^8)g(q^6,q^{9})-q\cdot g(q^4,q^{11})g(q^3,q^{12})-
	q^2g(q,q^{14})g(q^{3},q^{12})\nonumber\\
	&&-\:q\cdot g(q^6,q^9)g(q^{2},q^{13})+q^2g(q^{5},q^{10})g(1,q^{15}).\label{phisum}
\end{eqnarray}

The right hand side of \eqref{phisum} is equal to the expansion of $\det(A)$ plus an additional term $q^2g(q^{5},q^{10})g(1,q^{15})$, which is identically $0$. To see this, we expand the product:
\begin{equation}
g(q^{5},q^{10})g(1,q^{15})=\sum_{(m,n)\in \mathbb{Z}^2}(-1)^{m+n} q^{m(15m+5)/2+n(15n+15)/2}.\label{zeroterm}
\end{equation}
Next, we choose the index 2 sublattice $L_1:=\text{span}_{\mathbb{Z}}\{(1,1),(1,-1)\},$ of $\mathbb{Z}^2$, which has cosets $L_1, (1,0)+L_1$.  The sum in \eqref{zeroterm} becomes, after simplification:
\begin{eqnarray*}
g(q^5,q^{10})g(1,q^{15})
	&=&f(q^{5},q^{25})f(q^{10},q^{20})-q^{10}f(q^{-10},q^{40})f(q^5,q^{25})\\
	&=&f(q^{5},q^{25})f(q^{10},q^{20})-f(q^{10},q^{20})f(q^{5},q^{25})\\
	&=&0.
\end{eqnarray*}
Therefore $D_3=\varphi(q)^2$.
\end{proof}

 Hence for the case $n=3$ the generating functions for the outer multiplicities simplify to:
\begin{eqnarray}
B_1(q)&=&\frac{g(q^7,q^8)-q\:g(q^2,q^{13})}{\varphi(q)},\label{n3}\\
B_2(q)&=&\frac{g(q^4,q^{11})+q\:g(q,q^{14})}{\varphi(q)}.\label{n4}
\end{eqnarray}

For $k \in \mathbb{Z}_{\ge 0}$, let
\begin{eqnarray*}
a(k)&:=&|\mathscr{C}_{3,2\Lambda_0-k\delta}|\\
	&=&\# \text{ of partitions } (\lambda_1^{f_1},\lambda_2^{f_2},\dots, \lambda_j^{f_j}) \text{ of }3k\text{ such that } 0<f_i<3\text{ for all }i\leq j,\\
	&&f_1\equiv \lambda_1 \pmod{3}, \text{ and } f_i+f_{i+1}+\lambda_{i}-\lambda_{i+1}\equiv0 \pmod{3} \text{ for all } i<j,\\
b(k)&:=&|\mathscr{C}_{3,2\Lambda_1-k\delta}|\\
	&=&\# \text{ of partitions } (\lambda_1^{f_1},\lambda_2^{f_2},\dots, \lambda_j^{f_j}) \text{ of }3k-2\text{ such that } 0<f_i<3\text{ for all }i\leq j,\\
	&&f_1\equiv \lambda_1 \pmod{3}, \text{ and } f_i+f_{i+1}+\lambda_{i}-\lambda_{i+1}\equiv0 \pmod{3} \text{ for all } i<j,\\
c(k)&:=&\# \text{ of partitions of }k \text{ into parts } \not \equiv 0, \pm7 \pmod{15} - \\
	&&\# \text{ of partitions of }k-1 \text{ into parts } \not \equiv 0, \pm 2 \pmod{15},\\
d(k)&:=&\# \text{ of partitions of }k-1 \text{ into parts } \not \equiv 0, \pm4 \pmod{15} + \\
	&&\# \text{ of partitions of }k-2 \text{ into parts } \not \equiv 0, \pm 1 \pmod{15}.
\end{eqnarray*}
Then, by Theorem 3.1 and equations \eqref{n3} and  \eqref{n4} we obtain the following partition identities, which we believe to be new:
\begin{theorem}
$a(k)=c(k)$, for $k\geq 0$, and $b(k)=d(k)$, for $k\geq 1$.
\end{theorem}


\begin{thebibliography}{99}
\bibitem[A]{A}Andrews, George E.  \emph{The Theory of Partitions},  Cambridge University Press.  1984.
\bibitem[F]{F}Feingold, Alex J.  Tensor products of certain modules for the generalized Cartan matrix Lie algebra $A_{1}^{(1)}$,  \emph{Comm. Algebra}  \textbf{9} (1981) no. 12, 1323--1341.
\bibitem[H]{H}Humphreys, James. \emph{Introduction to Lie Algebras and Representation Theory},  Springer-Verlag.  1972.
\bibitem[HK]{HK}Hong, J.  Kang, S.-J.  \emph{Introduction to Quantum Groups and Crystal Bases},  American Mathematical Society.  2002.
\bibitem[JMMO]{JMMO}Jimbo, M. Misra, K. C.  Miwa, T.  Okado, M.  Combinatorics of Representations of $U_q(\widehat{\mathfrak{sl}}(n))$ at $q=0$,  \emph{Commun. Math. Phys.} \textbf{136} (1991), 543--556.
\bibitem[K1]{K1}Kac, V.  Infinite Dimensional Lie Algebras and Dedekind's $\eta$-function,  \emph{Functional Anal. Appl.} \textbf{8} (1974), 68--70.
\bibitem[K2]{K}Kac, V.  \emph{Infinite Dimensional Lie Algebras, $3^{\text{rd}}$ edition},  Cambridge University Press.  1990.
\bibitem[Ka1]{Ka1}Kashiwara, M. Crystalizing the $q$-analogue of universal enveloping algebras, \emph{Commun. Math. Phys.} \textbf{133} (1990), 249-260.
\bibitem[Ka2]{Ka2}Kashiwara, M. On crystal bases of the $q$-analogue of universal enveloping algebras, \emph{Duke Math. J.} \textbf{63} (1991), 465-516.
\bibitem[Li]{Li}Littelmann, P. A Littlewood-Richardson rule for symmetrizable Kac-Moody algebras, \emph{Invent. Math.} \textbf{116} (1994), 329-346.
\bibitem[Lu]{Lu}Lusztig, G. Canonical bases arising from quantized enveloping algebras, \emph{J. Amer. Math. Soc.} \textbf{3} (1990), 447-498.
\bibitem[M]{M} Misra, K. C.  Level Two Standard $\widetilde{A}_{n}$ modules,  \emph{J. Algebra} \textbf{137} (1991), 56--76.
\bibitem[MM]{MM} Misra, K. C., Miwa, T.   Crystal bases for the basic representation of $U_q(\widehat{\mathfrak{sl}}(n))$,  \emph{Commun. Math. Phys.} \textbf{134} (1990), 79-88.
\bibitem[MSZ]{MSZ}McLaughlin, J.  Sills, A. V.  Zimmer, P.  Rogers-Ramanujan-Slater Type Identities,  \emph{Electronic J. Combinatorics} DS15 (2008), 1--59,\\  (http://www.combinatorics.org/Surveys/ds15.pdf ).
\bibitem[OSS]{OSS} Okado, M. Schilling, A. Shimozono, M.   A tensor product theorem related to perfect crystals,  \emph{J. Algebra} \textbf{267} (2003), 212-245.
\end{thebibliography}
\end{document}